\documentclass[12pt]{amsart}
\usepackage{amsmath, amssymb}

\newcommand{\dbar}{\ensuremath{\overline\partial}}
\newcommand{\dbarstar}{\ensuremath{\overline\partial^*}}
\newcommand{\C}{\ensuremath{\mathbb{C}}}

\newcommand{\Z}{\ensuremath{\mathbb{Z}}}

\makeatletter     
\newcommand{\sumprime}{\if@display\sideset{}{'}\sum%
            \else\sum'\fi}    
\makeatother

\newtheorem{theorem}{Theorem}
\newtheorem{proposition}[theorem]{Proposition}
\newtheorem{lemma}[theorem]{Lemma}

\DeclareMathOperator{\dom}{Dom}

\begin{document}

\title[]{Semi-classical analysis of Schr\"{o}dinger operators
and compactness in the $\dbar$-Neumann problem}

\author{Siqi Fu}
\address{Department of Mathematics\\
         University of Wyoming\\
         Laramie,  WY 82071}
\email{sfu@uwyo.edu}
\author{Emil J.~Straube}
\address{Department of Mathematics\\
         Texas A \& M University\\
         College Station,  TX 77843} 	
\email{straube@math.tamu.edu}
\date {December 27, 2001.}
\thanks{The authors are supported in part by NSF grants
DMS-0070697 and DMS-0100517 respectively.}


\begin{abstract} We study the asymptotic behavior, in a
``semi-classical limit'', of the first eigenvalues (i.e. the
groundstate energies) of a class of Schr\"{o}dinger operators with
magnetic fields and the relationship of this behavior with compactness
in the $\bar\partial$-Neumann problem on Hartogs domains
in $\C^2$. \end{abstract}

\maketitle

\section{Introduction}\label{Introduction}

Let $\Omega$ be a bounded pseudoconvex domain in $\C^{n}$. The complex
Laplacian $\square_{q}$ is the operator $\dbar\dbarstar +
\dbarstar\dbar$, acting as an (unbounded) self-adjoint operator on
$\mathcal{L}^{2}_{0,q)}(\Omega)$, the space of $(0,q)$-forms with
coefficients in $\mathcal{L}^{2}(\Omega)$. It is a classical result of
H\"{o}r\-mander \cite{Hormander65} that $\square_{q}$ has a bounded
inverse. This inverse is the $\dbar$-Neumann operator $N_{q}$. The
$\dbar$-Neumann operator is closely related to solving
the $\dbar$-equation and thus plays a central role in several complex
variables. It is also of considerable interest from the point of view
of partial differential equations, where it provides a prototype (of
the solution) of an elliptic problem with non-coercive boundary
conditions. For a detailed survey of the $\mathcal{L}^{2}$-Sobolev
regularity theory of the $\dbar$-Neumann problem, we refer the reader
to~\cite{BoasStraube99}. In particular, it is known that global
regularity holds in many cases, but not all~\cite{Christ96}.
A question closely related to global regularity is that of compactness
of the $\dbar$-Neumann operator. This question is of interest in its
own right for a number of reasons; see~\cite{FuStraube01} for a
discussion of various aspects of the problem. In the context of global
regularity, the relevance stems from a theorem of Kohn
and Nirenberg~\cite{KohnNiren65} which implies that if $N_{q}$ is
compact in $\mathcal{L}^{2}_{(0,q)}(\Omega)$, then it is globally
regular in the sense that it preserves the $\mathcal{L}^{2}$-Sobolev
spaces. Catlin~\cite{Catlin84} demonstrated that compactness provides
indeed a viable route to global regularity for the $\dbar$-Neumann
problem, the link being his concept of property (P). He
showed that property(P) implies compactness (hence global regularity),
and that it can be verified on large classes of
domains. More recently, compactness is also being studied as a
property not only stronger than global regularity, but one that is
more robust and less subtle, and hence should be more amenable to a
reasonable characterization in terms of properties of the boundary.

In this paper, we relate
property(P) and compactness of the $\dbar$-Neumann operator on
complete pseudoconvex Har\-togs domains in $\mathbb{C}^{2}$ to the
asymptotic behavior of the groundstate energy of certain families of
Schr\"{o}dinger operators. It is well known that $\dbar$- and
related problems on such domains can be studied by means of the
corresponding weighted problem on the base domain, see for example
\cite{Ligocka89}, \cite{Berndtsson94} and their references. In turn,
studying the $\dbar$-equation in weighted $\mathcal{L}^{2}$-spaces on
planar domains leads naturally to Schr\"{o}dinger operators on these
domains \cite{Christ91},\cite{Berndtsson96}. We show that compactness
of the $\dbar$-Neumann operator and property(P) on the Hartogs domain
are characterized by the asymptotic behavior, in a ``semi-classical
limit'', of the lowest eigenvalues (the groundstate energies) of
certain magnetic Schr\"{o}dinger operators on the base domain and
their non-magnetic counterparts, respectively.

To state the main result of this paper, we need to introduce some
notation and recall some terminology. A compact set $K\subset \C^n$
is said to satisfy property ($P$) if for
every positive number $M$, there exists a neighborhood $U$ of $K$
and a $C^2$-smooth function $\lambda$ on $U$, 
$0\le \lambda \le 1$, such that for all $z$ in $K$, the
smallest eigenvalue of the Hermitian  form $\left(\partial^2\lambda
(z)/\partial z_j\partial\bar z_k \right)^n_{j, k=1}$ is at least $M$.
Let $D$ be a bounded domain
in $\C$ and let $\phi\in C^2(\overline D)$.  Let 
$S_\phi =-[(\partial_x +i \phi_y)^2
+(\partial_y-i\phi_x)^2] + \Delta \phi$ be a magnetic Schr\"{o}dinger
operator and $S^0_\phi=-\Delta +\Delta\phi$ be
the corresponding non-magnetic Schr\"{o}dinger operator.  Let
$\lambda_{\phi}(D)$ and $\lambda^0_{\phi}(D)$ be the first
eigenvalues of the Dirichlet realization of $S_\phi$ and 
$S^0_\phi$ on $D$ respectively. (See Section 2 below for details.) We
will also use the notation $\lambda(D)$ for $\lambda_{\phi \equiv 0}$,
that is, for the lowest eigenvalue of minus the Dirichlet Laplacian.

\begin{theorem} \label{Theorem 1}
Let $\Omega=\{ (z, w)\in \C^2;\ \ z\in  D, \ |w|<e^{-\phi(z)} \}$
be a smooth bounded complete pseudoconvex Hartogs domain in
$\C^2$. Suppose that $b\Omega$ is strictly pseudoconvex on
$b\Omega\cap\{w=0\}$.  Then \begin{enumerate}
\item $b\Omega$ satisfies property ($P$) if and only
if ${\displaystyle \lim_{n\to\infty}}
\lambda^0_{n\phi}(D) = \infty$.
\item $N$ is compact if and only if
${\displaystyle \lim_{n\to\infty}}\lambda_{n\phi}(D)=\infty$.
\end{enumerate}
\end{theorem}

The necessity in Part 2 in Theorem 1 is implicit in
Proposition~2 of Matheos' paper~\cite{Matheos97}. Note that $S_{n\phi}
= -n^{2}[(\frac{1}{n}\partial_x +i
\phi_y)^2+(\frac{1}{n}\partial_y-i\phi_x)^2] + n\Delta \phi$. Letting
$n$ tend to infinity is thus analogous, in a sense, to letting
``Planck's constant'' $h=1/n$ tend to zero. Study of the latter situation
is often referred to as semi-classical analysis (see
e.g.~\cite{Helffer88}, chapter 1).

Global regularity is not an issue for the domains we study here: the
$\dbar$-Neumann problem is globally regular on any smooth bounded
complete pseudoconvex Hartogs domain in $\C^{2}$
(~\cite{BoasStraube89}).

Sibony (\cite{Sibony87}, see also~\cite{Sibony91}) undertook a
systematic study of property(P), under the name of B-regularity, on
arbitrary compact sets in $\C^{n}$. Some of this work is also
discussed in section 3 of~\cite{FuStraube01}. In particular,
in the situation of Theorem~1, $b\Omega$ satisfies
property(P) (in $\C^{2}$) if and only if $W:=\{z \in D|\Delta
\phi = 0 \}$ satisfies property(P) in the plane
(\cite{Sibony87}, p.310, see also section 5 below).

We refer the reader to~\cite{FuStraube01} for a detailed discussion
of compactness in the $\dbar$-Neumann problem. As mentioned above, on
a bounded pseudoconvex domain, property(P) implies compactness of the
$\dbar$-Neumann operator. For sufficiently smooth domains,
see~\cite{Catlin84}; the second author observed in~\cite{Straube97}
that no boundary regularity at all is needed. In light of Theorem 1, a
quantitative way to look at this result in the case of Hartogs domains
in $\C^{2}$ is through (a special case of) Kato's inequality (see
Proposition~\ref{KatoInequality} below): $\lambda_{n\phi}(D)\ge
\lambda^0_{n\phi}(D)$. It would be of considerable interest, both from
the point of view of the $\dbar$-Neumann problem and from that of
Schr\"{o}dinger operators, to determine whether or not conversely, the
(equivalent) properties in part (2) of Theorem~\ref{Theorem 1} imply
those in part(1). For an example of a continuous (but
nonsmooth) subharmonic $\phi$ where $\lim_{n\rightarrow \infty}
\lambda_{n\phi}(D) = \infty$, but $\lim_{n\rightarrow \infty}
\lambda_{n\phi}^{0}(D) < \infty$, see~\cite{ChristFu01}.

Recently, McNeal~\cite{McNeal01} showed that a variant of property(P)
still implies compactness of the $\dbar$-Neumann operator. For the
domains we consider here, this variant turns out to be equivalent to
property(P); we discuss this in the appendix (Section~5).

\bigskip
\noindent{\bf Acknowledgment.} Part of this work was done while the 
first author visited Princeton University on an AMS Centennial
Research Fellowship.  The author is indebted to Professors 
J.~J.~Kohn  and E.~Lieb for stimulating discussions. 
He also thanks the Mathematics Department for hospitality.

\section{Schr\"{o}dinger operators}\label{Schrodinger}

In this section, we collect some facts about the Schr\"{o}dinger
operators that arise in the setting of Theorem 1. The reader may
in addition consult~\cite{Berndtsson96}.

Let $D$ be a bounded
domain in $\C$ and let $\phi (z)\in C^2(\overline{D})$. Let $\bar
L_\phi = e^{-\phi}\frac{\partial}{\partial \bar z} (e^{\phi}\cdot)
=\partial_{\overline{z}} +\phi_{\bar z}$ be the first order
differential operator defined in the sense of distributions on
$\mathcal{L}^2(D)$ and let $L_\phi = -e^{\phi}\frac{\partial}{\partial
z} (e^{-\phi}\cdot)=-\partial_{z} +\phi_z$ be the (formal) adjoint of
$\bar L_\phi$.  The domain of the actual adjoint of $\bar L_\phi$  is
the Sobolev space $W^1_0(D)$. Note that $\overline{L}_{\phi}$ is just
$\partial/\partial \overline{z}$ conjugated by multiplication by
$e^{\phi}$.

Consider the closed, positive semi-definite sesquilinear form
$$
Q_\phi(u, v)=4(L_\phi u, L_\phi v)
$$ 
defined on $W^1_0(D)\times
W^1_0(D)\subset \mathcal{L}^2(D)\times \mathcal{L}^2(D)$.  Let $S_\phi$ be
the unique non-negative, self-adjoint, densely defined operator on
$\mathcal{L}^2(D)$ corresponding to $Q_\phi(u, v)$. For the connection
between quadratic forms and (unbounded) self-adjoint operators, see
for example~\cite{ReedSimon80}, section VIII.6. Then
$\dom(S_{\phi})=W^{1}_{0}(D) \cap W^{2}(D)$, and on this domain
\begin{equation}\label{MagneticSchrodinger} S_\phi =4 \bar L_\phi
L_\phi =-[(\partial_x +i \phi_y)^2 +(\partial_y-i\phi_x)^2] + \Delta
\phi . \end{equation}
This is the Dirichlet realization of the Schr\"{o}dinger operator on $D$
with magnetic potential $A=(-\phi_y, \phi_x)=-\phi_y dx + \phi_x dy$,
magnetic field $dA=\Delta\phi dx\wedge dy$,
and electric potential $V=\Delta\phi$. Since
$\dom(S_{\phi})=W^{1}_{0}(D) \cap W^{2}(D)$ embeds compactly into
$\mathcal{L}^{2}(D)$, $S_{\phi}$ has compact resolvent. By
construction, $S_{\phi}$ is the restriction to its domain of an
isomorphism of $W^{1}_{0}(D)$ onto $W^{-1}(D)$. Consequently, as an
(unbounded) operator on $\mathcal{L}^{2}(D)$, it is injective and
onto, and so has a bounded inverse (which moreover is then compact).
Let $S^0_\phi=-\Delta +\Delta\phi$ be the Schr\"{o}dinger operator
without magnetic potential corresponding to $S_{\phi}$ (with the same
domain as $S_{\phi}$); it too has compact resolvent. Because both
$S_{\phi}$ and $S^{0}_{\phi}$ have compact resolvent, the spectrum in
each case consists of a sequence of eigenvalues tending to infinity.
Let $\lambda_\phi(D)$ and $\lambda^0_\phi(D)$ be the first eigenvalues
of $S_\phi$ and $S^0_\phi$, respectively. Note that
\begin{align}\label{Equivalence} 
\lambda_\phi (D)&=\inf\left\{4\int_D
|L_\phi u|^2 dA/\int_D |u|^2 dA, \ \ u\in C^\infty_0(D),
u\not\equiv 0\right\} \notag\\
&=\inf\left\{4\int_D |u_z|^2 e^{2\phi} dA/ \int_D |u|^2 e^{2\phi} dA, \ \
u\in C^\infty_0(D), u\not\equiv 0\right\} .
\end{align}

The case of most interest to us is that of a subharmonic $\phi$
(i.e.$\Delta\phi \geq 0$); this corresponds to pseudoconvexity of the
Hartogs domain $\Omega$. It turns out that subharmonicity of $\phi$
can be characterized in terms of the behavior of the lowest
eigenvalues $\lambda^{0}_{t\phi}(D)$ and $\lambda_{t\phi}(D)$. The
equivalence of (1) and (3) below is in~\cite{Berndtsson93}.

\begin{proposition}\label{Subharmonic}
The following are equivalent:
\begin{enumerate} 
\item $\Delta \phi \ge 0$;
\item $\liminf_{t\to \infty} \lambda_{t\phi}^{0} (D) >0$;
\item $\limsup_{t\to\infty} \lambda_{t\phi}(D) >0$.
\end{enumerate}
\end{proposition}

\begin{proof} (1) implies (2) because $S^{0}_{t\phi} \geq -\Delta$ (as
operators) when $\Delta\phi \geq 0$, and $-\Delta > 0$ (as operators).
(2) implies (3) because $\lambda_{t\phi}(D) \geq
\lambda^{0}_{t\phi}(D)$; for this, see
Proposition~\ref{KatoInequality} below. Finally, the proof that (3)
implies (1) can be found in~\cite{Berndtsson93}, proof of Proposition
1.5 and the remark immediately following that proof (bottom of page
212). \end{proof}

\medskip

\noindent{\bf Remark 1.}
1)Note that $\lambda(D)$ provides a lower bound on
$\liminf_{t\to \infty} \lambda_{t\phi}^{0}(D)$ that is independent of
$\phi$, for $\phi$'s that satisfy one of the properties in
Proposition~\ref{Subharmonic}. In particular, the quantities in
(2) and (3) in the proposition cannot be arbitrarily small positive
for $D$ given. $\lambda(D)$ itself can be estimated from below by
$\pi/|D|$, where $|D|$ denotes the area of $D$; this is a consequence
of the Poincar\'{e} inequality (see e.g. \cite{GilbargTrudinger98},
inequality (7.44)).

2)In general, one does not have $\limsup_{t\to\infty}
\lambda_{t\phi}(D)=\liminf_{t\to\infty} \lambda_{t\phi} (D)$.  For example,
if $D$ is an annulus and $\Delta \phi =0$  on $D$ then $\lambda_{t\phi}(D)$
is a periodic function of $t$ with minimum $\lambda(D)$ and maximum
$\lambda (D\setminus L)$ where $L$ is a path connecting the  
components of $bD$ (see~\cite{Helfferetc99}, Theorem 1.~1 and Remark
1.5~(vi)).

\medskip

The first part of the following proposition is a special
case of an inequality of Kato (\cite{Kato72}; see
also~\cite{Berndtsson96}). The second part goes back
to~\cite{LavineOCarroll77}. A detailed proof of the proposition is
in~\cite{Helffer88}, Lemma 7.2.1.2 and Theorem 7.2.1.1, to where we
refer the reader.

\begin{proposition}\label{KatoInequality}
$\lambda_\phi (D)\ge \lambda^0_\phi (D)$. Furthermore, equality
holds if and only if (1) $\Delta \phi=0$ on $D$ and (2)
$(1/2\pi)\int_\gamma A\in \Z$ for any simple closed smooth
curve $\gamma$ in $D$. In this case, $S_\phi (D)$ and
$S^0_\phi (D)$ are unitarily equivalent via $u
\leftrightarrow e^{ih}u$ for a $(\mathbb{R} \bmod 2\pi)$-valued $h$
with $dh=A$.
 \end{proposition}

The proof in \cite{Helffer88} uses the following
identity from \cite{LavineOCarroll77}, which can be obtained by
integration by parts:
\begin{equation}\label{Integral2}
4\|L_\phi u\|^2 -\lambda^0_\phi (D) \| u\|^2
=\|(\partial_x +i\phi_y -u^0_x/u^0)u\|^2
+ \|(\partial_y - i\phi_x -u^0_y/u^0)u\|^2 .
\end{equation}

Here, $u_{0}$ denotes the (normalized) eigenfunction of $S_{\phi}^{0}$
to the eigenvalue $\lambda_{\phi}^{0}$. (This eigenvalue is known to
be simple, and $u_{0}$ is known to be zero free in $D$, see
e.g.~\cite{ReedSimon80}, also the discussion in 7.2.1 in
~\cite{Helffer88}.) In particular, the difference $\lambda_{\phi}(D) -
\lambda_{\phi}^{0}(D)$ is the infimum of the right hand side of
\eqref{Integral2} over all $u \in C^{\infty}_{0}(D)$ with $||u||=1$.
Replacing $\phi$ by $n\phi$ in \eqref{Integral2} gives a semi-explicit
expression for $\lambda_{n\phi}(D) - \lambda_{n\phi}^{0}(D)$, which is
the quantity of concern in the question of whether or not property(P)
and compactness of the $\dbar$-Neumann operator are actually
equivalent for the domains in $\C^{2}$ considered in
Theorem~\ref{Theorem 1}. We call the expression semi-explicit because
it involves taking the infimum and it involves the eigenfunction
$u_{0}$, which depends on $n$. This latter dependence may possibly be
mitigated  by passing to subsequences which converge in appropriate
senses; compare the discussion in Section~3. Nonetheless,
formula~\eqref{Integral2} notwithstanding, it appears that determining
whether or not $\lambda_{n\phi}(D)$ tending to infinity implies that
$\lambda_{n\phi}^{0}(D)$ tends to infinity is a nontrivial matter.
Note that if $W:=\{z \in D| \Delta\phi(z)=0 \}$ has non-empty
interior, both eigenvalues are bounded above, for each
$n$, by the corresponding eigenvalues on, say, a disc
contained in the interior of $W$. On such a disc, the magnetic and
non-magnetic eigenvalues agree, by Proposition~\ref{KatoInequality},
and moreover, the non-magnetic eigenvalue is the same as that of
$-\Delta$. So both sequences are bounded above (and thus fail to
converge to infinity). Alternatively, if $W$ contains a disc,
$b\Omega$ contains an analytic disc, and both property(P) and
compactness of the $\dbar$-Neumann operator fail
(see~\cite{FuStraube01} for details). 
In contrast, when $W$ has empty \emph{fine}
interior (see Section~3), then $b\Omega$ satisfies property(P) and the
$\dbar$-Neumann operator on $\Omega$ is compact so that both sequences
of eigenvalues tend to infinity (in view of Theorem~\ref{Theorem 1}).
Consequently, the case of interest is that of $W$ with empty Euclidean
interior, but non-empty fine interior. (When $W$ has non-empty fine
interior, there exists \emph{some} smooth
subharmonic function $\psi$, such that $W=\{z\in\C|\Delta\psi =
0\}$ and $\lim_{n\rightarrow \infty}\lambda_{n\psi}(D) < \infty$, see
Remark 3 below.)

We point out that when no smoothness
restrictions are placed on the boundary of $\Omega$, compactness in
the $\dbar$-Neumann problem does not imply property(P). The domain
$\{(z,w)\in \C^{2} | 0<|z|<1, |z|^{2}+|w|^{2} < 1 \}$, obtained by
deleting from the unit ball the variety $\{z=0\}$, has an analytic
disc in its boundary (\emph{a fortiori}, the boundary does not have
property(P)), yet its $\dbar$-Neumann operator is compact
(see~\cite{FuStraube01}, example on page 150 preceding Proposition
4.1). The point is that the $\mathcal{L}^{2}$-theory does not detect
the deletion of the variety $\{z=0 \}$, and as a result, the
$\dbar$-Neumann operator inherits compactness from the $\dbar$-Neumann
operator on the unit ball. Recently, Christ and Fu
(~\cite{ChristFu01}) have constructed an example of a continuous
subharmonic $\phi$ with $\bigtriangledown\phi \in \mathcal{L}^{2}(D)$,
$\Delta\phi \in \mathcal{L}^{1}(D)$, and $\lim_{n \rightarrow \infty}
\lambda_{n\phi}(D) = \infty$, but $\lim_{n \rightarrow \infty}
\lambda_{n\phi}^{0}(D) < \infty$.

That magnetic Schr\"{o}dinger operators majorize
their non-magnetic counterparts in some appropriate sense, such as
Kato's inequality, is generally referred to as diamagnetism, and the
opposite direction (usually in terms of more general so called Pauli
operators) is called paramagnetism. The property in question in the
previous paragraph may thus be viewed as a paramagnetic property of
the family of Schr\"{o}dinger operators $\{S_{n\phi}|n\in\mathbb{N}
\}$ and their non-magnetic counterparts $\{S_{n\phi}^{0}|n\in\mathbb{N}
\}$. It appears that what is known in the theory of
Schr\"{o}dinger operators in this direction concerns cases that, when
specialized to the context of Theorem~\ref{Theorem 1}, cover
situations that are well understood from the point of view of the
$\dbar$-Neumann problem. For example, when the magnetic field is
constant a result of Lieb (which is in terms of more general Pauli
operators, see~\cite{AvronHerbstSimon78}) implies that
$\lambda_{\phi}(D) \leq \lambda_{2\phi}^{0}(D)$. Lieb's result was
later generalized by Avron and Seiler to the case where the magnetic
fields are given by certain polynomials (~\cite{AvronSeiler79}). The
ideas in the proofs of these results actually work when
$\phi(z)=\sum_{j=1}^{m}|h_{j}(z)|^{2}$.

\begin{proposition}\label{Paramagnetism}
If $\phi=\sum_{j=1}^m |h_j(z)|^2$ where $h_j(z)$ are holomorphic
on $D$, then
\begin{equation}\label{ParaIneq}
\lambda_\phi (D) \le \lambda^0_{2\phi}(D) .
\end{equation}
\end{proposition}
\begin{proof} Let $g$ be a real-valued eigenfunction of $S^0_{2\phi}$
corresponding to $\lambda^0_{2\phi}=\lambda^0_{2\phi}(D)$.  
For $\zeta\in \C^m$, we let
$H(z, \zeta)=-\sum_{j=1}^m (h_j(z)\zeta_j + |\zeta_j|^2)$,
$\Psi(z, \zeta)=e^{-\phi+H(z, \zeta)}$, and $f=g\Psi$.  It follows
that $S_\phi (f)=\lambda^0_{2\phi} f + 4 (2\phi_z - H_z)\Psi g_{\bar z}$.
Therefore,
\begin{align*}
(S_\phi(f), f)&=\lambda^0_{2\phi} \|f\|^2 -2 \int_D 
\frac{\partial |\Psi|^2}{\partial z}\frac{\partial g^2}{\partial \bar z}\\
        &=\lambda^0_{2\phi} \|f\|^2 +2 \int_D \frac{\partial^2 |\Psi|^2}{\partial
z\partial \bar z} g^2\\
        &=\lambda^0_{2\phi} \|f\|^2 + 2\int_D |f|^2 
                            (|2\phi_z -H_z|^2 -2\phi_{z\bar z}) \\
	&=\lambda^0_{2\phi}\|f\|^2 +2 \sum_{j, k=1}^m 
                             \int_D h_{jz}\bar h_{kz}
		\frac{\partial^2 |\Psi|^2}{\partial \bar \zeta_j\partial \zeta_k} g^2 .
\end{align*}
Denote $I(\zeta)$ the last term above.  It follows from the divergence theorem
that $\int_{\C^m} I(\zeta) =0$.  Therefore, there is a $\zeta_0\in
\C^m$ such that $I(\zeta_0)\le 0$.
We then conclude the proof.
\end{proof}

\noindent{\bf Remark 2.}
Proposition~\ref{Paramagnetism} does not
hold for all $\phi$ with $\Delta\phi\ge 0$.  For example, if
$D=\{z\in\C; 1/2 <|z|<2\}$ and $\Delta\phi=0$ on $D$, then
$\lambda^0_{2\phi}(D)=\lambda (D)=\lambda_{\phi}^{0}(D)$, which is
strictly less than $\lambda_\phi (D)$ unless $(1/2\pi)\int_{|z|=1}
A\in \Z$ (Proposition~\ref{KatoInequality}). By a limiting process,
one can in fact find examples such that $D$ is the unit disc and
$\Delta \phi \ge 0$ on $D$ but \eqref{ParaIneq} fails.

\section{Some potential theory}

Recall that the fine topology on $\C$ is the weakest topology
so that every subharmonic function is continuous. A general reference
for the basic facts about the fine topology in $\C$ is~\cite{Helms69}.
We use $int_{f}$ to denote the interior in the fine topology.

The Dirichlet problem for (minus) the Laplacian can be formulated on
finely open sets; see~\cite{Fuglede99}, section 3, and the
references there for this formulation. The resulting theory
inherits many features of the classical theory, but avoids some of its
problems having to do with ``stability'' of sets (see~\cite{Fuglede99},
Remark 2.4). What matters for us is the behavior of the first
eigenvalue under a decreasing sequence of finely
open sets. If $U$ is finely open, we still use $\lambda(U)$ to denote
this eigenvalue. Then, if $\{U_{j}\}_{j=0}^{\infty}$ is a decreasing
sequence of bounded finely open sets
in $\C$, $\{\lambda(U_{j})\}_{j=0}^{\infty}$ is
increasing (as in the classical case), and (unlike the
classical case) $\lim_{j \rightarrow \infty}\lambda(U_{j}) =
\lambda(int_{f}(\cap_{j}U_{j}))$ (\cite{Fuglede99}, Theorem 2, part
$1^{\circ}$).

The next proposition combines work of Fuglede and Sibony. We need it in
the proof of part (1) of Theorem~\ref{Theorem 1}.

\begin{proposition}\label{FineInterior}
Let $K$ be a compact subset of $\C$. The following are equivalent:
\begin{enumerate}
\item $K$ satisfies property ($P$).
\item $K$ has empty fine interior.
\item $K$ supports no non-zero function in $W^{1}_{0}(\C)$.
\item For any open sets $U_j$ such that
$K\subset\subset U_{j+1}\subset\subset U_j$ and
$\cap_{j=1}^\infty \overline{U}_j =K$, $\lambda
(U_j)\to\infty$ as $j\to\infty$. \end{enumerate} \end{proposition}

\begin{proof}
The equivalence of (1) and (2) is~\cite{Sibony87}, Proposition 1.11.
Let $\{U_{j}\}_{j=0}^{\infty}$ be a sequence of (Euclidean) open sets
with $K \subset \subset U_{j+1}\subset \subset U_{j}$ and
$\cap_{j}\overline{U_{j}} = K$. If (2) holds, then $\{0\} =
W^{1}_{0}(int_{f}K) = \cap_{j=1}^\infty W^{1}_{0}(U_{j})$
(~\cite{Fuglede99}, Lemma 1.1,(ii)), which gives (3). If (3) holds,
(4) must hold. If not,there would exist a sequence of functions
$\{u_{j}\}_{j=1}^{\infty}$, $u_{j} \in W^{1}_{0}(U_{j})$, with
$\| u_{j} \| = 1$ and $\| \nabla u_{j}\| \leq const.$ (use
\eqref{Equivalence} for $\phi \equiv 0$), for a suitable sequence
$\{U_{j}\}_{j=1}^{\infty}$. Passing to a subsequence that converges
both weakly in $W^{1}_{0}(\C)$ and in norm in $\mathcal{L}^{2}$ of a
neighborhood of $K$ would yield a non-zero element of $W^{1}_{0}(\C)$
that is supported on $K$, contradicting (3). Finally, (4) implies (2)
because $\lim_{j \rightarrow \infty} \lambda(U_{j}) =
\lambda(int_{f}(\cap_{j}U_{j})) = \lambda(int_{f}K)$, and the last
quantity is finite if $int_{f}K \neq \emptyset$, see~\cite{Fuglede99},
Theorem 2, part $1^{\circ}$.

\end{proof}

\noindent{\bf Remark 3.}
There are also characterizations of sets with empty fine interior in
terms of logarithmic capacity (~\cite{Helms69}, chapter 10, section 5)
and in terms of Brownian motion (\cite{Doob84}, section 2.IX.15). Our
work shows that such a characterization can also be given in terms of
non-magnetic Schr\"{o}dinger operators: a compact set $K \subseteq \C$
has empty fine interior if and only if for every smooth subharmonic
function $\phi$ on a domain $D$ with $K\subset\subset D$, such that
$K \supseteq \{z\in \C|\Delta\phi=0\}$, $\lim_{n \rightarrow
\infty}\lambda_{n\phi}^{0}(D) = \infty$. If $K$ has empty fine
interior, then combining Proposition~\ref{FineInterior}
and (the proof of) part (1) of Theorem~\ref{Theorem 1} shows that
$\lim_{n \rightarrow \infty}\lambda_{n\phi}^{0}(D) = \infty$ for the
$\phi$'s under consideration. Conversely, the authors have shown
(~\cite{FuStraube01}, Theorem 4.2) that if $K$ has non-empty fine
interior, there exists a smooth, bounded, pseudoconvex, complete
Hartogs domain in $\C^{2}$ whose weakly pseudoconvex boundary points
project onto $K$ and whose $\dbar$-Neumann operator is not compact.
Moreover, the Hartogs domain can be chosen to satisfy the assumptions
in Theorem~\ref{Theorem 1}. Consequently, the resulting function
$\phi$ on the base domain $D$ satisfies $\lim_{n \rightarrow
\infty}\lambda_{n\phi}(D) < \infty$ (Theorem~\ref{Theorem 1}). \emph{A
fortiori}, $\lim_{n \rightarrow\infty}\lambda_{n\phi}^{0}(D) < \infty$
(Proposition~\ref{KatoInequality}).

\section{Proof of Theorem 1}

Let $W=\{z \in D| \Delta \phi(z)= 0 \}$. Because $b\Omega$ is
strictly pseudoconvex near the points of the boundary where $w=0$, $W$ 
is a compact subset of $D$. To prove part (1) of Theorem 1, we use
that $b\Omega$ satisfies property(P) if and only if $W$ satisfies
property(P) as a set in $\C$ (~\cite{Sibony87}, p. 310, see also
section 4 below). In turn, $W$ satisfies property(P) if and only it
satisfies (4) in Proposition~\ref{FineInterior}. Fix a sequence
$\{W_{j}\}_{j=1}^{\infty}$ of open subsets of $D$ such that $W=
\cap_{j}\overline{W_{j}}$, and $W\subset\subset W_{j+1}\subset\subset
W_{j}\subset\subset D$.

Assume $W$ satisfies (4) in Proposition~\ref{FineInterior}. Let
$\eta_j\in C^\infty_0(W_j)$, $0\le \eta_j\le 1$, and $\eta_j=1$ on
$W_{j+1}$. For any $u\in C^\infty_0(D)$, $j \in \mathbb{N}$,
\begin{align}
(S^0_{n\phi}u, u) &=\|\nabla u\|^2 + n\|\sqrt{\Delta \phi}
u\|^2 \notag\\
&\ge \frac1{2} \|\nabla
(u\eta_j)\|^2 + \frac{n}2 \|\sqrt{\Delta\phi}u\|^2 \notag\\
&\ge
\frac1{2}\lambda(W_j)\|u\eta_j\|^2 + \frac{n}2
\|\sqrt{\Delta\phi} u\|^2\notag \\
&\ge \frac{1}{2}\lambda(W_j)\|u\|^2\notag \end{align}
when $n$ is sufficiently large. By assumption,
$\lambda(W_{j}) \rightarrow \infty$ as $j \rightarrow \infty$.
Therefore $\lim_{n\to\infty}\lambda^0_{n\phi}(D)=\infty$. This
finishes one direction in the proof of part (1) of
Theorem~\ref{Theorem 1}.

For the other direction, observe that for all
$(n,j) \in \mathbb{N} \times \mathbb{N}$, $\lambda_{n\phi}^{0}(D) \leq
\lambda_{n\phi}^{0}(W_{j})$ (by the monotonicity with respect to the
domain of the eigenvalue). Also, for $u$ $\in C^{\infty}_{0}(W_{j})$,
\[ (S_{n\phi}^{0}u,u)=(- \Delta u+n\Delta\phi u,u) \leq (- \Delta
u,u)+(u,u) \]
if $j$ is big enough relative to $n$ so that $|n\Delta\phi| \leq 1$ on
$W_{j}$. Consequently, $\lambda(W_{j}) \geq
\lambda_{n\phi}^{0}(W_{j})-1 \geq \lambda_{n\phi}^{0}(D)-1$ if $j$ is
big enough relative to $n$. It follows that $lim_{j \rightarrow
\infty} \lambda(W_{j})$ 
$ = \infty$, since $\lim_{n \rightarrow
\infty}\lambda_{n\phi}^{0}(D) = \infty$. Since the sequence $W_{j}$
is arbitrary, this concludes the proof of part (1) of
Theorem~\ref{Theorem 1}.

We now prove the necessity in Part (2) of Theorem 1.  As noted before,
this also follows from Proposition 2 in \cite{Matheos97}.  We provide
a proof that does not require any regularity of $b\Omega$. (In this
case, $\lambda_{n\phi} (D)$ is defined by the second equality
in~\eqref{Equivalence}.)  We use the fact that compactness of $N$ is
equivalent to compactness of Kohn's canonical solution operator
$S=\dbarstar N$, which is in turn equivalent to the following
compactness estimates: For any $\epsilon>0$, there exists
$C_\epsilon>0$ such that \begin{equation}\label{CES}
\|Su\|^2\le \epsilon \|u\|^2 + C_\epsilon \|u\|^2_{-1}
\end{equation}
for all $u\in L^2_{(0, 1)}(\Omega)$( see~\cite{FuStraube01}, Lemma
1.1).

Let $\beta\in C^\infty_0(D)$ and let $u_n=\beta(z)w^n d\bar z$ and 
$f_n(z, w)=S(u_n)$.  Then $f_n(z, w)=g_n(z)w^n$ and $\partial
g_{n}(z)/\partial \bar z =\beta(z)$.  Plugging this into~\eqref{CES}
and using the fact that $\|g_n(z)w^n\|^2_{-1, \Omega}
\le (1/n^2) \|g_n(z)\|^2_\Omega$, we obtain that there exists $N_\epsilon>0$ such that
$\|g_n(z)w^n\|^2\le \epsilon \|\beta(z) w^n d\bar z\|^2$ when $n>N_\epsilon$.
Therefore, 
$$
\int_{D} |g_n(z)|^2 e^{-2(n+1)\phi(z)} dA(z) \le \epsilon \int_{D}
|\beta(z)|^2e^{-2(n+1)\phi(z)} dA(z).
$$
Duality gives for $u\in C^{\infty}_{0}(D)$

\begin{align}\notag
\int_{D} |u(z)|^2 e^{2(n+1)\phi(z)} &=\sup\{ |\langle u,
\beta\rangle|^2; \ \beta\in C^\infty_0 (D), \ \int_{D} |\beta|^2
e^{-2(n+1)\phi(z)} \le 1 \}\\ &\le \sup\{|\langle u_z,
g_n\rangle|^2; \ \int_{D} |g_n|^2 e^{-2(n+1)\phi} \le
\varepsilon\} \notag \\ &\le \varepsilon\int_{D} |u_z|^2
e^{2(n+1)\phi} . \notag
\end{align}
The middle inequality follows from consideration of the special
$g_{n}$ associated in the previous paragraph to a $\beta \in
C^\infty_0 (D)$. In view of~\eqref{Equivalence}, this concludes the
proof of necessity.

We now proof the sufficiency. Since compactness of the $\dbar$-Neumann
operator is a local property (see~\cite{FuStraube01}, Lemma 1.2; the
direction we need here follows from a simple partition of unity
argument) and since by assumption $b\Omega$ is strictly pseudoconvex
in a neighborhood of $b\Omega\cap \{w=0\}$, we need  only establish
compactness estimates (~\cite{FuStraube01}, Lemma 1.1) for forms whose
support is away from $b\Omega\cap\{w=0\}$. Moreover, by the interior
elliptic regularity of $\dbar \oplus \dbarstar$, we need only consider
forms whose support is close to $b\Omega$. Choose $\hat
D\subset\subset D$ such that $b\Omega$ is strictly pseudoconvex on a
neighborhood of the part of $b\Omega$ over the outside of $\hat D$.

We work for a moment on $\hat D \times S^{1}$ ($S^{1}$ is the unit
circle). Denoting the variables on $\hat D \times S^{1}$ by $(z,t)$,
let $L=\partial_z +i\phi_z \partial_t$. We use $||| \cdot |||$ to
denote norms on $\hat D \times S^{1}$. We will prove that for every
$\epsilon >0$, there exists $C_\epsilon>0$ such that
\begin{equation}\label{CEL}
|||u|||^2\le \epsilon (|||Lu|||^2 +|||\bar Lu|||^2) + C_\epsilon
|||u|||^2_{-1}
\end{equation}
for  $u\in C^\infty_0 (\hat D\times
S^{1})$. By the assumption on the eigenvalues $\lambda_{n\phi}(D)$,
there exists $N_\epsilon>0$ such that when $n>N_\epsilon$, $$
\|v\|^2\le \epsilon \|L_{n\phi} v\|^2, \quad \text{for all}\ v\in
C^\infty_0(\hat D). $$ (Note that $\lambda_{n\phi}(\hat D) \geq
\lambda_{n\phi}(D)$.) Taking conjugates, we obtain that when
$n<-N_\epsilon$, $$ \|v\|^2\le \epsilon \|\bar L_{n\phi} v\|^2, \quad
\text{for all}\ v\in C^\infty_0(\hat D). $$ Therefore, when
$|n|>N_\epsilon$  $$ \|v\|^2\le \epsilon (\|L_{n\phi} v\|^2+\|\bar
L_{n\phi} v\|^2), \quad \text{for all}\ v\in C^\infty_0(\hat D). $$
For $u\in C^\infty_0 (\hat D\times S^{1})$, write
$$
u=\sum_{n=-\infty}^\infty u_n(z) e^{int}
$$ where
$u_n(z)=(1/2\pi)\int_0^{2\pi} u(z, e^{it}) e^{-int} dt\in
C^\infty_0 (\hat D)$.  Then \[ Lu =
\sum_{n=-\infty}^{\infty}(-L_{n\phi}u_{n})e^{int}, \] and
\begin{align}
\frac{1}{2\pi}|||u|||^2 &=\frac{1}{2\pi}\sum_{n=-\infty}^\infty
|||u_n|||^2 = \sum_{n=-\infty}^\infty \|u_n\|^2
\notag \\         &\le \sum_{|n|\le N_\epsilon} \|u_n\|^2 +
\epsilon \sum_{|n|>N_\epsilon} (\|L_{n\phi} u_n \|^2 +\|\bar L_{n\phi}
u_n\|^2 ) \notag \\         &=\frac{\epsilon}{2\pi}(|||Lu|||^2+|||\bar
L u|||^2) + \sum_{|n|\le N_\epsilon} \left(\|u_n\|^2
-\epsilon(\|L_{n\phi} u_n\|^2 + \|\bar L_{n\phi} u_n\|^2)\right)
.\notag \end{align}
The last sum in the above inequalities is less
than or equal to $$
C_\epsilon\sum_{|n|\le N_\epsilon} \|u_n\|^2_{-1}
$$
for some sufficiently large $C_\epsilon$, depending only on
$\epsilon$. This is because $\forall n, L_{n\phi}$ and
$\overline{L}_{n\phi}$ have a compact inverse (see section 2), which
implies $||u_{n}||^{2} \leq \epsilon(||L_{n\phi}u_{n}||^{2} +
||\overline{L}_{n\phi}u_{n}||^{2}) +
C_{\epsilon}||u_{n}||_{-1}^{2}$ for a constant $C_{\epsilon}$. (This
is analogous to Lemma 1.1 in~\cite{KohnNiren65}, see also Lemma
1.1 in~\cite{FuStraube01}.) $C_{\epsilon}$ depends on $n$, but
because we are now only concerned with $n$'s satisfying $|n|\le
N_\epsilon$, $C_{\epsilon}$ may be chosen
depending only on $\epsilon$. The desired inequality \eqref{CEL} now
follows from the fact that the last sum above is controlled by
$|||u|||^2_{-1}$.

We now return to the setting of the Hartogs domain in
Theorem~\ref{Theorem 1}. For the part of the boundary over $\hat D$,
we may use as defining function the function $\rho(z,w) =
\frac{1}{2}\log(w\overline{w}e^{2\phi})$. For, say, $0 < r < 1$, the
level sets $M_{r}=\{\rho = -r\}$ are the surfaces $\{|w|^{2}=e^{2\phi
-2r}\}$. For $r$ fixed,
we use coordinates $(z,t)$ on $M_{r}$ via $(z,t) \leftrightarrow
(z,e^{-\phi(z)-r+it})$. Denote by $L_{1}$ the usual complex tangential
field of type (1,0) given by $\rho_{z}\partial_{w} -
\rho_{w}\partial_{z}$. A computation shows that when restricted to
$M_{r}$, $2wL_{1}$ becomes $\partial_{z} +i\phi_{z}\partial_{t}$,
which is the operator $L$ considered in the previous paragraph. Let
now $u$ be a smooth function supported above $\hat D$ and sufficiently
close to $b\Omega$. Denote by $d\sigma_{r}$ the surface measure on
$M_{r}$. Using that $dV$ in $\C^{2}$ is comparable to $d\sigma_{r}
\,dr$ (on $supp(u)$), and $d\sigma_{r}$ is comparable to $dV(z) \,dt$,
uniformly in $r$, we obtain from~\eqref{CEL} \begin{align*}
||u||^{2}=\int_{\Omega}|u|^{2} &\simeq
\int_{0}^{1}(\int_{M_{r}}|u|^{2}d\sigma_{r})dr \\
&\lesssim
\epsilon\int_{0}^{1}(\int_{M_{r}}(|Lu|^{2} +
|\overline{L}u|^{2})d\sigma_{r})dr +
C_{\epsilon}\int_{0}^{1}||u||^{2}_{-1,M_{r}}dr \\
&\lesssim
\epsilon\int_{0}^{1}(\int_{M_{r}}(|L_{1}u|^{2}+|\overline{L}_{1}u|^{2}
)d\sigma_{r})dr + C_{\epsilon}||u||_{-1}^{2} \\
&\lesssim
\epsilon(||L_{1}u||^{2}+||\overline{L}_{1}u||^{2})+C_{\epsilon}||u||_{
-1}^{2}. 
\end{align*}
Here, as usual, $\lesssim$ indicates ``less than or equal to, up to a constant
factor that is independent of $\epsilon$''. Let now $\alpha =
a_{1}d\overline{z}+a_{2}d\overline{w} \in
C^{\infty}_{(0,1)}(\overline{\Omega}) \cap \dom\dbarstar$, with
support above $\hat D$ and close to $b\Omega$. Changing $C_{\epsilon}$ if
necessary, we get from the estimate above

\[ ||\alpha||^{2} \leq
\epsilon(||L_{1}\alpha||^{2}+||\overline{L}_{1}\alpha||^{2}) +
C_{\epsilon}||\alpha||_{-1}^{2}, \]
where $L_{1}$ and $\overline{L_{1}}$ act componentwise on forms,
as usual.

We next invoke maximal estimates (~\cite{Derridj78}, Th\'{e}or\`{e}me
3.1): in $\C^{2}$,
$||L_{1}\alpha||^{2}+||\overline{L}_{1}\alpha||^{2}$ is controlled by
$||\dbar\alpha||^{2}+||\dbarstar\alpha||^{2}$. (Actually, the
statement in~\cite{Derridj78} includes the term $||\alpha||^{2}$, but
this term is now well known to be bounded by
$||\dbar\alpha||^{2}+||\dbarstar\alpha||^{2}$; alternatively, we may
absorb it into the left hand side.) The result of applying the maximal
estimates is (again, $C_{\epsilon}$ may have to be increased):
\[ ||\alpha||^{2} \leq
\epsilon(||\dbar\alpha||^{2}+||\dbarstar\alpha||^{2}) +
C_{\epsilon}||\alpha||_{-1}^{2}. \]
This is the required compactness estimate. The proof of
Theorem~\ref{Theorem 1} is complete.
\bigskip

\noindent{\bf Remark 4.}
The assumption in Theorem~\ref{Theorem 1} that $\Omega$ is strictly
pseudoconvex near the boundary of the base is not essential. It
suffices for example that the boundary is of finite type
(~\cite{D'Angelo93}) near points of $b\Omega \cap \{w=0\}$. One can
then replace $W$ by the (Euclidean) closure of $int_{f}(W)$ in the
above proofs (compare~\cite{Fuglede99}). This set will be relatively
compact in $D$ because $\Delta\phi$ vanishes to infinite order at fine
interior points of $W$ (see~\cite{Helms69}, Corollary 10.5 or Theorem
10.14). We leave the details to the reader.

\section{Appendix}

In this section, we show that on the domains considered in
Theorem~\ref{Theorem 1}, property($P$) and property($\tilde P$) are
actually equivalent.

We first recall the definition of property ($\tilde P$) by
McNeal in \cite{McNeal01}.  A compact set $K$
in $\C^n$ is said to satisfy property ($\tilde P$)
if for any $M>0$, there exists a neighborhood 
$U$ of $K$ and $g\in C^2(U)$ such that
\begin{itemize}
\item[(1)] $|\langle\partial g, X\rangle|^2\le L_g(X)$;
\item[(2)] $L_g(X)\ge M |X|^2$.
\end{itemize}
Here $\langle\cdot, \cdot\rangle$ is the pairing between
a form and a vector and 
$L_g(X)=\dbar\partial g(X, \overline{X})$. (1) is equivalent
to $-e^{-g}$ being plurisubharmonic in $U$ (see the discussion
in~\cite{McNeal01}).

\begin{lemma}\label{Sibony}
Let $\Omega$ be a smooth bounded complete
pseudoconvex Hartogs domain in $\C^2$.  Assume
that $b\Omega$ is strictly pseudoconvex at 
the base.  Then $b\Omega$ satisfies property ($P$)
if and only if it satisfies property ($\tilde P$).
\end{lemma}

\begin{proof}  It is easy to see that property ($P$) always 
implies property ($\tilde P$)(~\cite{McNeal01}): if $\lambda$ is the
function in the definition of property($P$),  it suffices (modulo
a normalization) to consider the function $g=e^\lambda$. The other
direction follows by combining Lemma~\ref{base1} and Lemma~\ref{base2}
below. \end{proof}

Let $\Omega=\{ (z, w); \ \ z\in D, |w|<e^{-\phi}\}$. Then $\Delta\phi\ge 0$
and the weakly pseudoconvex points correspond to the set of base
points $W=\{z\in D|\Delta\phi=0\}$ for. Note that
$W\subset\subset D$.

\begin{lemma}\label{base1}
Let $K$ be a compact subset of $\C$.  Then $K$ satisfies property ($P$)
if and only if it satisfies property ($\tilde P$).
\end{lemma}
\begin{proof}
We only have to show that property ($\tilde P$)
implies property ($P$).  In light of Proposition~\ref{FineInterior},
it suffices to show that for any open sets $U_j$ such that 
$K\subset\subset U_{j+1}\subset\subset U_j$ and $\cap_{j=1}^\infty
\overline{U}_j
=K$, $\lim_{j\to\infty}\lambda(U_j)=\infty$.

For any $M>0$, there exists a neighborhood $U$ of 
$K$ and $g\in C^2(U)$ such that $|g_z|^2\le g_{z\bar z}$
and $g_{z\bar z}\ge M$ on $U$.  Assume that $j_0$ is
sufficiently large so that $U_{j_0}\subset\subset U$.
It follows from an easy integration by parts that
\[
\int_{U_{j_0}} |u_z -\tfrac12 g_z u|^2 dA
=\tfrac12 \int_{U_{j_0}} g_{z\bar z} |u|^2 dA
+\int_{U_{j_0}} |u_{\bar z} +\tfrac12 g_{\bar z} u|^2 dA 
\]
for any $u\in C^\infty_0 (U_{j_0})$.   The left hand side
of the above equation is bounded from above by 
$3\|u_z\|^2 +\tfrac38\|g_z u\|^2$ while the right hand
side is bounded from below by $\tfrac12\int g_{z\bar z} |u|^2 dA$.
Therefore,
\[
\int_{U_{j_0}} |u_z|^2 dA \ge \tfrac1{24}\int_{U_{j_0}} g_{z\bar z} |u|^2 dA
\ge \frac{M}{24}\int_{U_{j_0}} |u|^2 dA .
\]
Hence $\lambda(U_j)\ge \lambda(U_{j_0})\ge M/6$ when $j\ge j_0$.
This concludes the proof of Lemma~\ref{base1}.
\end{proof}

\begin{lemma}\label{base2}
Assumptions as in Lemma~\ref{Sibony}. Then
\begin{enumerate}
\item$b\Omega$ satisfies property ($P$)  if and only if $W$
satisfies property ($P$).
\item$b\Omega$ satisfies property ($\tilde P$)  if and only if $W$
satisfies property ($\tilde P$).
\end{enumerate}
\end{lemma}

\begin{proof}
Part (1) may be found in~\cite{Sibony87}, page 310.

To prove (2), first note that if $W$ satisfies property($\tilde P$),
it satisfies property($P$) (Lemma~\ref{base1}), hence so does
$b\Omega$, by part(1). But then $b\Omega$ also satisfies
property($\tilde P$), by the discussion above.

The proof of the other direction is completely analogous to the proof
of the corresponding direction in (1). We are indebted to Nessim
Sibony for a private communication (~\cite{Sibony95}) on the details
of the argument in~\cite{Sibony87}. Fix $M >0$. Let $g$ be the
corresponding plurisubharmonic function from the definition of
property($\tilde P$). Replacing $g$
by$(1/2\pi)\int_{0}^{2\pi}g(z,we^{i\theta})d\theta$, we may assume
that $g$ is invariant under rotations in the $w$ variable. Consider
$h(z):= g(z, e^{-\phi(z)})$, defined in a neighborhood of
$\overline{D}$. Then, for a sufficiently small neighborhood $U$ of $W$

\begin{equation}\label{PtildeH}
h_{z\overline{z}} \geq M \; ;\; |h_{z}|^{2} \leq h_{z\overline{z}}.
\end{equation}
This is a matter of computation. This computation can be somewhat
simplified by first observing that the function $g_{1}(z,w):=g(z,
e^{w})$, defined in a neighborhood of the set $\{(z,w)\in\C^{2}|
z\in D, w+ \overline{w} = -2\phi(z)\}$, also satisfies (1) and (2) in
the definition of property ($\tilde P$), with $M$ replaced by, say,
$\tilde M = (\min\{e^{-2|\phi(z)|-1} | z\in U\})M$, where $U$ is a
suitable neighborhood of $W$ (after shrinking the neighborhood where
$g_{1}$ is defined, independently of $M$). Now $h(z)=g_{1}(z,
-\phi(z))$; also note that since $g$ is invariant under rotations in
the $w$ variable, $g_{1}$ is independent of the imaginary part of $w$,
that is, $(g_{1})_{w} \equiv (g_{1})_{\overline{w}}$.
It follows that
\[ h_{z} = (g_{1})_{z} -2(g_{1})_{w}\phi_{z} = \langle\partial g_{1},
X\rangle \]
and
\[ h_{z\overline{z}} = L_{g_{1}}(X)-2
(g_{1})_{w}\phi_{z\overline{z}}, \]
where $X=(1,-2\phi_{z})$. Consequently,
~\eqref{PtildeH} is satisfied at points of $W$ (where
$\phi_{z\overline{z}} = 0$), up to replacing $M$ by $\tilde M$.
Rescaling $h$ (for example, replacing $h$ by $h/2$) allows one to
conclude~\eqref{PtildeH} for $z$ in a small enough neighborhood of $W$
(by continuity). \end{proof}

\bigskip

\bibliography{survey}

\begin{thebibliography}{XXXXX}


%
%

\bibitem[AHS78]{AvronHerbstSimon78}
J.~Avron, I.~Herbst, and B.~Simon,
\emph{Schr\"{o}dinger operators with magnetic fields I. General interactions},
Duke Math. Journal \textbf{45} (1978), 847-883.

\bibitem[AS79]{AvronSeiler79}
J.~E.~Avron and R.~Seiler,
\emph{Paramagnetism for nonrelativistic electrons and {E}uclidean massless
{D}irac particles}, Physical Review Letters, \textbf{42} (1979), 931-934.

\bibitem[Be93]{Berndtsson93}
Bo Berndtsson,
\emph{A smooth pseudoconvex domain in $\C^2$ for which $L^\infty$-estimates
for $\dbar$ do not hold},
Ark. Mat. \textbf{31} (1993), 209-218.

\bibitem[Be94]{Berndtsson94}
\bysame
\emph{Some recent results on estimates for the $\dbar$-equation},
Contributions to Complex Analysis and Analytic Geometry
(H.~Skoda and J.~M.~Tr\'{e}preau, eds.), Aspects of Mathematics,
vol.~E26, Vieweg, 1994, pp. 27-42.

\bibitem[Be96]{Berndtsson96}
\bysame
\emph{\dbar \ and {S}chr\"{o}dinger operators},
Math. Z. \textbf{221} (1996), 401-413.

\bibitem[BS89]{BoasStraube89}
Harold~P. Boas and Emil~J. Straube,
\emph{Complete Hartogs domains in $\C^{2}$ have regular Bergman and
Szeg\"{o} projections},
Math. Z. \textbf{201} (1989), 441-454.

\bibitem[BS99]{BoasStraube99}
\bysame
\emph{Global regularity of the $\dbar$-Neumann problem:
a survey of the $L^2$-{S}obolev theory}, Several Complex
Variables (M. Schneider and Y.-T. Siu, eds.), MSRI Publications,
vol.~37, Cambridge University Press, 1999, pp. 79-111.


\bibitem[Ca84]{Catlin84}
David Catlin, 
\emph{Global regularity of the {$\overline\partial$}-{Neumann}
problem}, Complex Analysis of Several Variables (Yum-Tong Siu, ed.),
Proceedings of Symposia in Pure Mathematics, no.~41, American Mathematical
Society, 1984, pp.~39--49.


\bibitem[Chr91]{Christ91}
Michael Christ,
\emph{On the $\dbar$-equation in weighted $L^2$-norms in $\C$},
J. of Geom. Analysis \textbf{1}(1991), 193-230.

\bibitem[Chr96]{Christ96}
\bysame
\emph{Global $C^{\infty}$ irregularity of the $\dbar$-Neumann problem
for worm domains}, J. Amer. Math. Soc.\textbf{9}, Nr.4 (1996),
1171-1185.

\bibitem[ChrF01]{ChristFu01}
Michael Christ and Siqi Fu,
in preparation.

\bibitem[CFKS87]{Cycon87}
H.~L.~Cycon, R.~G.~Froese, W.~Kirsch, and B.~Simon,
\emph{Schr\"{o}dinger Operators}, with applications to quantum
mechanics, Springer, 1987.

\bibitem[DA93]{D'Angelo93}
John P.~D'Angelo,
\emph{Several Complex Variables and the Geometry of Real
Hypersurfaces},  Studies in Advanced Mathematics, CRC Press, 1993.


\bibitem[D78]{Derridj78}
M.~Derridj,
\emph{Regularit\'{e} pour $\dbar$ dans quelques domaines faiblement
pseudo-convexes}, J~.Diff.~Geometry \textbf{13}, Nr.4 (1978), 559-576.

\bibitem[Do84]{Doob84}
J.~L.~Doob,
\emph{Classical Potential Theory and Its Probabilistic Counterpart},
Grundlehren der mathematischen Wissenschaften \textbf{262}, Springer,
1984.

\bibitem[FS98]{FuStraube98}
Siqi Fu and Emil J.~Straube,
\emph{Compactness of the $\dbar$-{N}eumann problem on convex domains},
J. Func. Analysis. \textbf{159} (1998), 629-641.

\bibitem[FS01]{FuStraube01}
\bysame,
\emph{Compactness in the $\dbar$-{N}eumann problem}, Complex Analysis
and Geometry (J.~McNeal, ed.), Ohio State Math. Res. Inst. Publ.
\textbf{9} (2001), 141-160.

\bibitem[F72]{Fuglede72}
Bent Fuglede, 
\emph{Finely Harmonic Functions}, Lecture Notes in Math. 
vol.~289, Springer, 1972.

\bibitem[F99]{Fuglede99}
\bysame, 
\emph{The {D}irichlet {L}aplacian on finely open sets}, 
Potential Anal. \textbf{10} (1999), 91--101.

\bibitem[GT98]{GilbargTrudinger98}
D.~Gilbarg and N.~S.~Trudinger,
\emph{Elliptic Partial Differential Equations of Second Order}, second
edition,
Grundlehren der mathematischen Wissenschaften  \textbf{224}, Springer,
1998.


\bibitem[Hel88]{Helffer88}
B.~Helffer,
\emph{Semi-Classical Analysis for the Schr\"{o}dinger Operator and
Applications}, Lecture Notes in Math. No.~1336, Springer, 1988.

\bibitem[HHHO99]{Helfferetc99}
B.~Helffer, M.~Hoffmann-Ostenhof, T.~Hoffmann-Ostenhof, and
M.~P.~Owen,
\emph{Nodal sets for groundstates of Schr\"{o}dinger operators with
zero magnetic field in non simply connected domains},
Commun. Math. Phys. \textbf{202} (1999), 629-649.


\bibitem[He69]{Helms69}
L.~L.~Helms,
\emph{Introduction to Potential Theory},
Wiley-Interscience, 1969.


\bibitem[H65]{Hormander65}
Lars H{\"{o}}rmander,
\emph{{$L\sp{2}$} estimates and existence theorems for
  the {$\overline\partial$} operator}, Acta Mathematica \textbf{113} (1965),
  89--152.

\bibitem[Ka72]{Kato72}
T.~Kato,
\emph{Schr\"{o}dinger operators with singular potentials}, Israel J. Math. \textbf{13}
(1972), 135-148.

\bibitem[KN65]{KohnNiren65}
J.~J.~Kohn and L.~Nirenberg, \emph{Non-coercive boundary value
problems}, Commun. Pure and Applied Math. \textbf{18}
(1965), 443-492.


\bibitem[LO77]{LavineOCarroll77}
Richard Lavine and Michael O'Carroll,
\emph{Ground state properties and lower bounds for energy
levels of a particle in a uniform magnetic field and external
potential},
J. of Math. Physics \textbf{18} (1977), 1908-1912.

\bibitem[L89]{Ligocka89}
Ewa Ligocka,
\emph{On the Forelli-Rudin construction and weighted Bergman
projections},
Studia Math. \textbf{94} (1989), 257-272.

\bibitem[M97]{Matheos97}
Peter Matheos, \emph{A {H}artogs domain with no analytic discs in the 
boundary for which the $\dbar$-{N}eumann problem is not compact},
preprint, 1997 (to appear in J. of Geom. Analysis).

\bibitem[McN01]{McNeal01}
Jeffery D.~McNeal,
\emph{A sufficient condition for compactness of the $\dbar$-Neumann
operator},
preprint, 2001.

\bibitem[RS80]{ReedSimon80}
Michael Reed and Barry Simon,
\emph{Methods of Modern Mathematical Physics}, vol.1 (Functional
Analysis), Academic Press, 1980.



\bibitem[Si87]{Sibony87}
N.~Sibony,
\emph{Une classe de domaines pseudoconvexes}, Duke Mathematical
Journal \textbf{55} (1987), no.~2, 299--319.


\bibitem[Si91]{Sibony91}
\bysame,
\emph{Some aspects of weakly pseudoconvex domains},
Several Complex Variables and Complex Geometry, (E.~Bedford,
J.~P.~D'Angelo, R.~E.~Greene, and S.~G.~Krantz, editors), Proceedings
of Symposia in Pure Mathematics, vol.52, part 1, American Mathematical
Society, 1991, pp.199-231.

\bibitem[Si95]{Sibony95}
\bysame,
private correspondence.

\bibitem[St97]{Straube97}
Emil J.~Straube,
\emph{Plurisubharmonic functions and subellipticity of the $\dbar$-{N}eumann
problem on non-smooth domains}, Math. Res. Lett. \textbf{4}~(1997), 459-467.



\end{thebibliography}

\providecommand{\bysame}{\leavevmode\hbox to3em{\hrulefill}\thinspace}

\end{document}